\documentclass[11pt,reqno]{amsart}

\usepackage{graphics}
\usepackage{amssymb, a4wide}
\usepackage{amsmath}
\usepackage{graphics}
\usepackage{graphicx}
\usepackage{epsfig}
\newtheorem{definition}{Definition}  

\newtheorem{example}{Example}

\title{Space-like Loxodromes on the Canal Surfaces in Minkowski 3-Space}

\author{NILGUN SONMEZ}
\address[N. Sonmez]{Afyon Kocatepe University,
Department of Mathematics,
03200, Afyonkarahisar, Turkey}
\email{nceylan@aku.edu.tr}

\author{Murat Babaarslan}
\address[M. Babaarslan]{Bozok University,
Department of Mathematics,
66100, Yozgat, Turkey}
\email{murat.babaarslan@bozok.edu.tr}



\date{\today}

\begin{document}

\parindent 0mm
\parskip 2mm

\maketitle

\begin{abstract}
In this paper, we obtain the differential equations of the space-like loxodromes on the non-degenerate canal surfaces depending on the causal characters of these canal surfaces and their meridians in Minkowski 3-space. Also we give an example by using Mathematica computer programme.
\end{abstract}

\section{\protect\small Introduction}

Loxodromes are special curves which cuts all meridians on the Earth's surface at a constant angle. Thus loxodromes are usually used in navigation. Noble \cite{noble} found the equations of the loxodromes on the rotational surfaces in Euclidean 3-space. Babaarslan and Munteanu \cite{babaarslan1} obtained the equations of time-like loxodromes on the rotational surfaces in Minkowski 3-space. After that the equations of space-like loxodromes on the rotational surfaces in Minkowski 3-space were given by Babaarslan and Yayli \cite{babaarslan2}. A canal surface in Euclidean 3-space is the envelope of a moving sphere whose trajectory of centers is a spine curve $\alpha(u)$ with varying radius $r(u)$. The analytic and algebraic properties of canal surfaces in Euclidean 3-space were given by Xu et al. \cite{xu}. A lot of object and structures can be represented by using canal surfaces, for examples; pipes, hoses, brass instruments, internal organs of the body in solid modeling, helical channel and tunnels. Cylinder, cone, torus, sphere, tubular surfaces and Dupin cyclide are some particular examples of canal surfaces. If we take the radius $r(u)$ as constant, then the canal surfaces reduce to tubular surfaces (\cite{babaarslan4},\cite{ucum}). Ucum and Ilarslan \cite{ucum} found the parametrizations of the canal surfaces for all Lorentz spheres which are pseudo sphere (De Sitter space) $\mathbb{S}^2_{1}$, pseudo-hyperbolic sphere (hyperbolic plane) $\mathbb{H}^2$ and light-like cone $\mathbb{C}$. The differential equations of loxodromes on canal surfaces in Euclidean 3-space were obtained by Babaarslan \cite{babaarslan4}. We know that helicoidal surfaces are a natural generalization of the rotational surfaces. Loxodromes on helicoidal surfaces in Euclidean 3-space were studied by Babaarslan and Yayli \cite{babaarslan3}. Differential equations of the space-like and time-like loxodromes on helicoidal surfaces in Minkowski 3-space were found by Babaarslan and Kayacik (\cite{babaarslan5}, \cite{babaarslan6}).

In this paper, we investigate the differential equations of space-like loxodromes on the non-degenerate canal surfaces in Minkowski 3-space which were obtained by Ucum and Ilarslan \cite{ucum}. Also we give an example by using Mathematica.

\section{\protect\small Preliminaries}

In this section we give necessary concepts related to curves and surfaces in Minkowski 3-space $\mathbb{E}^3_{1}$.

The Lorentzian scalar product of the vectors $u=(u_{1},u_{2},u_{3})$ and $v=(v_{1},v_{2},v_{3})$ in $\mathbb{E}^3_{1}$ is
\begin{eqnarray}
\label{eq1}
<u,v>=-u_{1}v_{1}+u_{2}v_{2}+u_{3}v_{3}.
\end{eqnarray}
Also the pseudo-norm of the vector $u\in\mathbb{E}^3_{1}$ is given by
\begin{eqnarray}
\label{eq2}
\|u\|=\sqrt{|\langle u,u\rangle|}.
\end{eqnarray}

An arbitrary vector $u\in\mathbb{E}^3_{1}$ is (or its causal character is)
\begin{enumerate} \item[ i.]space-like if $\langle u,u\rangle>0$ or $u=0$,\\ \item[ii.] time-like if $\langle u,u\rangle<0$,\\ \item[iii.] light-like (null) if $\langle u,u\rangle=0$ and $u\neq0$.
\end{enumerate}
Let $\alpha:I\rightarrow \mathbb{E}^3_{1}$ be a regular curve in $\mathbb{E}^3_{1}$, where $I\subset \mathbb{R}$ is an open interval. The curve $\alpha$ is called
\begin{enumerate} \item[ i.] space-like if $\langle \dot{\alpha},\dot{\alpha}\rangle>0$,\\ \item[ii.] time-like if $\langle \dot{\alpha},\dot{\alpha}\rangle<0$,\\ \item[iii.] light-like if $\langle \dot{\alpha},\dot{\alpha}\rangle=0$  (see \cite{lopez}).
\end{enumerate}

Let $S:U\rightarrow\mathbb{E}^3_{1}$ be a smooth immersed surface in $\mathbb{E}^3_{1}$, where $U\subset\mathbb{R}^2$ is an open set. $S$ is non-degenerate if its first fundamental form is non-degenerate. Then $S$ is (or its causal character is)
\begin{enumerate} \item[ i.] space-like if its first fundamental form is a Riemannian metric,\\ \item[ii.] time-like if its first fundamental form is a Lorentzian metric (see \cite{sipus}).
\end{enumerate}
For examples; the pseudo-hyperbolic sphere (hyperbolic plane)
\begin{eqnarray}
\label{eq7}
\mathbb{H}^2=\{p\in\mathbb{E}^3_{1}|<p,p>=-1\}
\end{eqnarray}
is a space-like surface.

Also, the pseudo sphere (De Sitter space)
\begin{eqnarray}
\label{eq7}
\mathbb{S}^2_{1}=\{p\in\mathbb{E}^3_{1}|<p,p>=1\}
\end{eqnarray}
is a time-like surface (see \cite{lopez}).

Let $\{S_{u},S_{v}\}$ be a local base of the tangent plane at each point of $S$. Then the first fundamental form of $S$ is
\begin{eqnarray}
\label{eq6}
I=ds^{2}=Edu^{2}+2Fdudv+Gdv^{2},
\end{eqnarray}
where $E=\langle S_{u},S_{u}\rangle$, $F=\langle S_{u},S_{v}\rangle$ and $G=\langle S_{v},S_{v}\rangle$ are the coefficients of first fundamental form of $S$.

By using these coefficients, we can give the causal characters of $S$. For example;

\begin{enumerate} \item[ i.] $S$ is a space-like  if and only if $\det(I)=EG-F^{2}>0$,\\ \item[ii.] $S$ is a time-like  if and only if $\det(I)=EG-F^{2}<0$ (see \cite{lopez}, \cite{ugail}).
\end{enumerate}

Also the arc-length of any curve on $S$ between $u_{1}$ and $u_{2}$ can be given by
\begin{eqnarray}
\label{eq7}
s=\bigg|\int_{u_{1}}^{u_{2}}\sqrt{\bigg|E+2F\frac{dv}{du}+G(\frac{dv}{du})^{2}\bigg|}du\bigg|
\end{eqnarray} (see \cite{babaarslan5}).

\section{Differential equations of the space-like loxodromes on the
space-like canal surfaces which have space-like meridians}

In this section, we find the differential equations of the space-like loxodromes on the space-like canal surfaces having space-like meridians (that is, $EG-F^{2}>0$ and $E>0$ for all $(u,v)$).

\begin{definition}
\label{d6} If $u$ and $v$ are space-like vectors in ${\mathbb{E}}_{1}^{3}$
which span a space-like plane. Then
\begin{equation}
\langle u,v\rangle =\Vert u\Vert \Vert v\Vert \cos \psi,\notag
\end{equation}%
where $\psi \in {\mathbb{R}}$ $(0\leq \psi \leq \pi )$ is the Lorentzian
space-like angle between $u$ and $v$ \cite{ratcliffe}.
\end{definition}

\begin{definition}
\label{d7} A space-like curve on a space-like canal surface in ${\mathbb{E}}%
^3_1$ which has space-like meridians is called a loxodrome if the curve cuts
all meridians at a constant Lorentzian space-like angle.
\end{definition}
\textbf{1.} Let us consider the following space-like canal surface $C$ which is given by
(3.1) in \cite{ucum}:
\begin{align}
C\left( u,v\right) =\alpha \left( u\right) +h\left( u\right) T\left(
u\right) +g\left( u\right) m_{1}\left( \sinh v\right) N\left( u\right)
+g\left( u\right) m_{2}\left( \cosh v\right) B\left( u\right),
\end{align}%
where $g\left( u\right) =r\left( u\right) \sqrt{1+\left( r^{\prime }\left(
u\right) \right) ^{2}}$,\ $h\left( u\right) =r\left( u\right) r^{\prime
}\left( u\right) $ and $m_{1},\ m_{2}\in \{-1,1\}$.

The coefficients of first fundamental form of the canal surface $C$ are
\begin{eqnarray*}
E &=&\langle C_{u},C_{u}\rangle =\left( k_{1}h\left( u\right)
-k_{2}m_{2}g\left( u\right) \cosh v+m_{1}g^{\prime }\left( u\right) \sinh
v\right) ^{2} \\
&&-\left( k_{2}m_{1}g\left( u\right) \sinh v-m_{2}g^{\prime }\left( u\right)
\cosh v\right) ^{2}+\left( 1-k_{1}m_{1}g\left( u\right) \sinh v+h^{\prime
}\left( u\right) \right)^{2}, \\
&&\text{ \ \ \ \ \ \ \ \ \ \ \ \ } \\
F &=&\langle C_{u},C_{v}\rangle =g\left( u\right) \left( k_{1}m_{1}h\left(
u\right) \cosh v-k_{2}m_{1}m_{2}g\left( u\right) +\left(
m_{1}^{2}-m_{2}^{2}\right) g^{\prime }\left( u\right) \cosh v\sinh v\right),
\\
&&\text{\ \ \ \ \ \ \ \ \ \ \ \ \ } \\
G &=&\langle C_{v},C_{v}\rangle =g^{2}\left( u\right) \left( m_{1}^{2}\cosh
^{2}v-m_{2}^{2}\sinh ^{2}v\right).
\end{eqnarray*}

By using these coefficients, the first fundamental form of $C$ is given by
\begin{eqnarray*}
ds^{2} &=&\{\left( k_{1}h\left( u\right) -k_{2}m_{2}g\left( u\right) \cosh
v+m_{1}g^{\prime }\left( u\right) \sinh v\right) ^{2} \\
&&-\left( k_{2}m_{1}g\left( u\right) \sinh v-m_{2}g^{\prime }\left( u\right)
\cosh v\right) ^{2}+\left( 1-k_{1}m_{1}g\left( u\right) \sinh v+h^{\prime
}\left( u\right) \right) ^{2}\}du^{2} \\
&& \\
&&+2\{g\left( u\right) \left( k_{1}m_{1}h\left( u\right) \cosh
v-k_{2}m_{1}m_{2}g\left( u\right) +\left( m_{1}^{2}-m_{2}^{2}\right)
g^{\prime }\left( u\right) \cosh v\sinh v\right) \}dudv \\
&& \\
&&+\left( g^{2}\left( u\right) \left( m_{1}^{2}\cosh ^{2}v-m_{2}^{2}\sinh
^{2}v\right) \right) dv^{2}.
\end{eqnarray*}
Let us assume that $\beta(t)$ is the image of a curve $(u(t),v(t))$ on the $(uv)$-plane under $C$. According to the local basis $\{C_{u},C_{v}\}$, the tangent vector $\beta^{\prime}(t)$ has the coordinates $(u^{\prime },v^{\prime })$ and the tangent vector $C_{u}$ has the coordinates $(1,0)$, Therefore, at the point $C(u,v)$, where the space-like loxodrome cuts the space-like meridians at a constant Lorentzian space-like angle, we get
\begin{eqnarray*}
\cos \psi &=&\dfrac{Edu+Fdv}{\sqrt{\Delta }} \\
&& \\
&=&\dfrac{\{\left( k_{1}h\left( u\right) -k_{2}m_{2}g\left( u\right) \cosh
v+m_{1}g^{\prime }\left( u\right) \sinh v\right) ^{2}}{\sqrt{\Delta }} \\
&&\dfrac{-\left( k_{2}m_{1}g\left( u\right) \sinh v-m_{2}g^{\prime }\left(
u\right) \cosh v\right) ^{2}+\left( 1-k_{1}m_{1}g\left( u\right) \sinh
v+h^{\prime }\left( u\right) \right) ^{2}\}}{\sqrt{\Delta }}du \\
&& \\
&&+\dfrac{\{g\left( u\right) \left( k_{1}m_{1}h\left( u\right) \cosh
v-k_{2}m_{1}m_{2}g\left( u\right) +\left( m_{1}^{2}-m_{2}^{2}\right)
g^{\prime }\left( u\right) \cosh v\sinh v\right) \}}{\sqrt{\Delta }}dv,
\end{eqnarray*}
where
\begin{eqnarray*}
\Delta &=&E^{2}du+2EFdudv+EGdv^{2} \\
&& \\
&=&\{\left( k_{1}h\left( u\right) -k_{2}m_{2}g\left( u\right) \cosh
v+m_{1}g^{\prime }\left( u\right) \sinh v\right) ^{2} \\
&&-\left( k_{2}m_{1}g\left( u\right) \sinh v-m_{2}g^{\prime }\left( u\right)
\cosh v\right) ^{2}+\left( 1-k_{1}m_{1}g\left( u\right) \sinh v+h^{\prime
}\left( u\right) \right) ^{2}\}^{2}du \\
&&+\{2g\left( u\right) \left( k_{1}m_{1}h\left( u\right) \cosh
v-k_{2}m_{1}m_{2}g\left( u\right) +\left( m_{1}^{2}-m_{2}^{2}\right)
g^{\prime }\left( u\right) \cosh v\sinh v\right) \\
&&\times\lbrack \left( k_{1}h\left( u\right) -k_{2}m_{2}g\left( u\right) \cosh
v+m_{1}g^{\prime }\left( u\right) \sinh v\right) ^{2}-\left(
k_{2}m_{1}g\left( u\right) \sinh v-m_{2}g^{\prime }\left( u\right) \cosh
v\right) ^{2} \\
&&+\left( 1-k_{1}m_{1}g\left( u\right) \sinh v+h^{\prime }\left( u\right)
\right) ^{2}]\}dudv \\
&&+\{g^{2}\left( u\right) \left( m_{1}^{2}\cosh ^{2}v-m_{2}^{2}\sinh
^{2}v\right) [\left( k_{1}h\left( u\right) -k_{2}m_{2}g\left( u\right) \cosh
v+m_{1}g^{\prime }\left( u\right) \sinh v\right) ^{2} \\
&&-\left( k_{2}m_{1}g\left( u\right) \sinh v-m_{2}g^{\prime }\left( u\right)
\cosh v\right) ^{2}+\left( 1-k_{1}m_{1}g\left( u\right) \sinh v+h^{\prime
}\left( u\right) \right) ^{2}]\}dv^{2}.
\end{eqnarray*}

From this equation, we obtain the following differential equation of the space-like loxodrome on the space-like canal surface having space-like meridians:
\begin{equation*}
\left( \cos ^{2}\psi EG-F^{2}\right) \left( \dfrac{dv}{du}\right) ^{2}-2\sin
^{2}\psi EF\dfrac{dv}{du}=\sin ^{2}\psi E^{2},
\end{equation*}%
that is
\begin{eqnarray*}
&&\{g^{2}\left( u\right) [-(k_{1}m_{1}h\left( u\right) \cosh
v-k_{2}m_{1}m_{2}g(u)+\left( m_{1}^{2}-m_{2}^{2}\right) g^{\prime }\left(
u\right) \cosh v\sinh v)^{2} \\
&&+\cos ^{2}\psi \left( m_{1}^{2}\cosh ^{2}v-m_{2}^{2}\sinh ^{2}v\right) ]%
\text{ }[\left( k_{1}h\left( u\right) -k_{2}m_{2}g(u)\cosh v+m_{1}g^{\prime
}\left( u\right) \sinh v\right) ^{2} \\
&&-\left( k_{2}m_{1}g\left( u\right) \sinh v-m_{2}g^{\prime }\left( u\right)
\cosh v\right) ^{2}+\left( 1-k_{1}m_{1}g\left( u\right) \sinh v+h^{\prime
}\left( u\right) \right) ^{2}]\}\left( \dfrac{dv}{du}\right) ^{2} \\
&& \\
&&-\{2\sin ^{2}\psi g\left( u\right) \left( k_{1}m_{1}h\left( u\right) \cosh
v-k_{2}m_{1}m_{2}g\left( u\right) +\left( m_{1}^{2}-m_{2}^{2}\right)
g^{\prime }\left( u\right) \cosh v\sinh v\right) \\
&&\times\lbrack \left( k_{1}h\left( u\right) -k_{2}m_{2}g\left( u\right) \cosh
v+m_{1}g^{\prime }\left( u\right) \sinh v\right) ^{2}-\left(
k_{2}m_{1}g\left( u\right) \sinh v-m_{2}g^{\prime }\left( u\right) \cosh
v\right) ^{2} \\
&&+\left( 1-k_{1}m_{1}g\left( u\right) \sinh v+h^{\prime }\left( u\right)
\right) ^{2}]\}\dfrac{dv}{du} \\
&& \\
&=&\sin ^{2}\psi \{\left( k_{1}h\left( u\right) -k_{2}m_{2}g\left( u\right)
\cosh v+m_{1}g^{\prime }\left( u\right) \sinh v\right) ^{2} \\
&&-\left( k_{2}m_{1}g\left( u\right) \sinh v-m_{2}g^{\prime }\left( u\right)
\cosh v\right) ^{2}+\left( 1-k_{1}m_{1}g\left( u\right) \sinh v+h^{\prime
}\left( u\right) \right) ^{2}\}^{2}.
\end{eqnarray*}

\textbf{2.} Let us consider the following space-like canal surface $C$
which is given by (3.2) in \cite{ucum}:
\begin{equation}
C\left( u,v\right) =\alpha \left( u\right) +h\left( u\right) T\left(
u\right) +g\left( u\right) m_{1}\left( \cosh v\right) N\left( u\right)
+g\left( u\right) m_{2}\left( \sinh v\right) B\left( u\right),
\end{equation}
where $%
\begin{array}{c}
g\left( u\right) =r\left( u\right) \sqrt{1+\left( r^{\prime }\left( u\right)
\right) ^{2}}%
\end{array}%
,%
\begin{array}{c}
h\left( u\right) =r\left( u\right) r^{\prime }\left( u\right)%
\end{array}%
$ and $m_{1},\ m_{2}\in \{-1,1\}$.

The coefficients of first fundamental form of the canal surface $C$ are
\begin{eqnarray*}
E &=&-\left( -k_{1}h\left( u\right)
+k_{2}m_{2}g\left( u\right) \sinh v+m_{1}g^{\prime }\left( u\right) \cosh
v\right) ^{2} \\
&&+\left( k_{2}m_{1}g\left( u\right) \cosh v+m_{2}g^{\prime }\left( u\right)
\sinh v\right) ^{2}+\left( 1-k_{1}m_{1}g\left( u\right) \cosh v+h^{\prime
}\left( u\right) \right) ^{2}, \\
&&\text{ \ \ \ \ \ \ \ \ \ \ \ \ } \\
F &=&g\left( u\right) \left( k_{1}m_{1}h\left(
u\right) \sinh v+k_{2}m_{1}m_{2}g\left( u\right) +\left(
-m_{1}^{2}+m_{2}^{2}\right) g^{\prime }\left( u\right) \cosh v\sinh v\right),
\\
&&\text{\ \ \ \ \ \ \ \ \ \ \ \ \ } \\
G &=&g^{2}(u)\left(-m_{1}^{2}\sinh^{2}v+m_{2}^{2}\cosh^{2}v\right).
\end{eqnarray*}
By using these coefficients, the first fundamental form of $C$ is
\begin{eqnarray*}
ds^{2} &=&\{-\left( -k_{1}h\left( u\right) +k_{2}m_{2}g\left( u\right) \sinh
v+m_{1}g^{\prime }\left( u\right) \cosh v\right) ^{2} \\
&&+\left( k_{2}m_{1}g\left( u\right) \cosh v+m_{2}g^{\prime }\left( u\right)
\sinh v\right) ^{2}+\left( 1-k_{1}m_{1}g\left( u\right) \cosh v+h^{\prime
}\left( u\right) \right) ^{2}\}du^{2} \\
&& \\
&&+2\{g\left( u\right) \left( k_{1}m_{1}h\left( u\right) \sinh
v+k_{2}m_{1}m_{2}g\left( u\right) +\left( -m_{1}^{2}+m_{2}^{2}\right)
g^{\prime }\left( u\right) \cosh v\sinh v\right) \}dudv \\
&& \\
&&+\{g^{2}(u)\left(-m_{1}^{2}\sinh^{2}v+m_{2}^{2}\cosh^{2}v\right)\}dv^{2}.
\end{eqnarray*}%

At the point $C(u,v)$, where the space-like loxodrome cuts the space-like meridians at a constant
Lorentzian space-like angle, we get

\begin{eqnarray*}
\cos \psi&=&\dfrac{\{-\left( -k_{1}h\left( u\right) +k_{2}m_{2}g\left( u\right) \sinh
v+m_{1}g^{\prime }\left( u\right) \cosh v\right) ^{2}}{\sqrt{\Delta }} \\
&&\dfrac{+\left( k_{2}m_{1}g\left( u\right) \cosh v+m_{2}g^{\prime }\left(
u\right) \sinh v\right) ^{2}+\left( 1-k_{1}m_{1}g\left( u\right) \cosh
v+h^{\prime }\left( u\right) \right) ^{2}\}}{\sqrt{\Delta }}du \\
&& \\
&&+\dfrac{\{g\left( u\right) \left( k_{1}m_{1}h\left( u\right) \sinh
v+k_{2}m_{1}m_{2}g\left( u\right) +\left( -m_{1}^{2}+m_{2}^{2}\right)
g^{\prime }\left( u\right) \cosh v\sinh v\right) \}}{\sqrt{\Delta }}dv,
\end{eqnarray*}
where
\begin{eqnarray*}
\Delta &=&\{-\left( -k_{1}h\left( u\right) +k_{2}m_{2}g\left( u\right) \sinh
v+m_{1}g^{\prime }\left( u\right) \cosh v\right) ^{2} \\
&&+\left( k_{2}m_{1}g\left( u\right) \cosh v+m_{2}g^{\prime }\left( u\right)
\sinh v\right) ^{2}+\left( 1-k_{1}m_{1}g\left( u\right) \cosh v+h^{\prime
}\left( u\right) \right) ^{2}\}^{2}du \\
&& \\
&&+\{2g\left( u\right) \left( k_{1}m_{1}h\left( u\right) \sinh
v+k_{2}m_{1}m_{2}g\left( u\right) +\left( -m_{1}^{2}+m_{2}^{2}\right)
g^{\prime }\left( u\right) \cosh v\sinh v\right)  \\
&&\times[-\left( -k_{1}h\left( u\right) +k_{2}m_{2}g\left( u\right) \sinh
v+m_{1}g^{\prime }\left( u\right) \cosh v\right) ^{2}+\left(
k_{2}m_{1}g\left( u\right) \cosh v+m_{2}g^{\prime }\left( u\right) \sinh
v\right) ^{2} \\
&&+\left( 1-k_{1}m_{1}g\left( u\right) \cosh v+h^{\prime }\left( u\right)
\right) ^{2}]\}dudv \\
&& \\
&&+\{\left(-m_{1}^{2}\sinh^{2}v+m_{2}^{2}\cosh^{2}v\right) g^{2}\left( u\right)[-\left( -k_{1}h\left( u\right) +k_{2}m_{2}g\left( u\right) \sinh
v+m_{1}g^{\prime }\left( u\right) \cosh v\right) ^{2} \\
&&+\left( k_{2}m_{1}g\left( u\right) \cosh v+m_{2}g^{\prime }\left( u\right)
\sinh v\right) ^{2}+\left( 1-k_{1}m_{1}g\left( u\right) \cosh v+h^{\prime
}\left( u\right) \right) ^{2}]\}dv^{2}.
\end{eqnarray*}

From this equation, we obtain the following differential equation of the space-like loxodrome:

\begin{eqnarray*}
&&\{-g^{2}\left( u\right) (k_{1}m_{1}h\left( u\right) \sinh
v+k_{2}m_{1}m_{2}g(u)+\left( -m_{1}^{2}+m_{2}^{2}\right) g^{\prime }\left(
u\right) \cosh v\sinh v)^{2} \\
&&+\cos ^{2}\psi \left(-m_{1}^{2}\sinh^{2}v+m_{2}^{2}\cosh^{2}v\right) g^{2}\left( u\right) [-(-k_{1}h\left( u\right) +k_{2}m_{2}g(u)\sinh
v+m_{1}g^{\prime }\left( u\right) \cosh v)^{2} \\
&&+\left( k_{2}m_{1}g\left( u\right) \cosh v+m_{2}g^{\prime }\left( u\right)
\sinh v\right) ^{2}+\left( 1-k_{1}m_{1}g\left( u\right) \cosh v+h^{\prime
}\left( u\right) \right) ^{2}]\}\left( \dfrac{dv}{du}\right) ^{2} \\
&& \\
&&-\{2\sin ^{2}\psi g\left( u\right) \left( k_{1}m_{1}h\left( u\right) \sinh
v+k_{2}m_{1}m_{2}g\left( u\right) +\left( -m_{1}^{2}+m_{2}^{2}\right)
g^{\prime }\left( u\right) \cosh v\sinh v\right) \\
&&\times[-\left( -k_{1}h\left( u\right) +k_{2}m_{2}g\left( u\right) \sinh
v+m_{1}g^{\prime }\left( u\right) \cosh v\right) ^{2}+\left(
k_{2}m_{1}g\left( u\right) \cosh v+m_{2}g^{\prime }\left( u\right) \sinh
v\right) ^{2} \\
&&+\left( 1-k_{1}m_{1}g\left( u\right) \cosh v+h^{\prime }\left( u\right)
\right) ^{2}]\}\dfrac{dv}{du} \\
&& \\
&=&\sin ^{2}\psi \{-\left( -k_{1}h\left( u\right) +k_{2}m_{2}g\left(
u\right) \sinh v+m_{1}g^{\prime }\left( u\right) \cosh v\right) ^{2} \\
&&+\left( k_{2}m_{1}g\left( u\right) \cosh v+m_{2}g^{\prime }\left( u\right)
\sinh v\right) ^{2}+\left( 1-k_{1}m_{1}g\left( u\right) \cosh v+h^{\prime
}\left( u\right) \right) ^{2}\}^{2}.
\end{eqnarray*}

\textbf{3.} Let us consider the following space-like canal surface $C$ which is given by (3.31) in \cite{ucum}:
\begin{equation}
C\left( u,v\right) =\alpha \left( u\right) -h\left( u\right) T\left(
u\right) +p\left( u\right) m_{1}\left( \cos v\right) N\left( u\right)
+p\left( u\right) m_{2}\left( \sin v\right) B\left( u\right),
\end{equation}%
where
$p\left( u\right) =r\left( u\right) \sqrt{\left( r^{\prime }\left(
u\right) \right) ^{2}-1}$,\ $h\left( u\right) =r\left( u\right) r^{\prime
}\left( u\right) $ and $m_{1},\ m_{2}\in \{-1,1\}$.

The coefficients of first fundamental form of the canal surface $C$ are
computed as
\begin{eqnarray*}
E &=&-\left( 1+k_{1}m_{1}p\left( u\right) \cos
v-h^{\prime }\left( u\right) \right) ^{2}+\left( k_{1}h\left( u\right)
+k_{2}m_{2}p\left( u\right) \sin v-m_{1}p^{\prime }\left( u\right) \cos
v\right) ^{2} \\
&&+\left( k_{2}m_{1}p\left( u\right) \cos v+m_{2}p^{\prime }\left( u\right)
\sin v\right) ^{2}, \\
&&\text{ \ \ \ \ \ \ \ \ \ \ \ \ } \\
F &=&p\left( u\right) \left( k_{1}m_{1}h\left(
u\right) \sin v+k_{2}m_{1}m_{2}p\left( u\right) +\left(
-m_{1}^{2}+m_{2}^{2}\right) p^{\prime }\left( u\right) \cos v\sin v\right), \\
&&\text{\ \ \ \ \ \ \ \ \ \ \ \ \ } \\
G &=&p^{2}\left( u\right) \left( m_{1}^{2}\sin
^{2}v+m_{2}^{2}\cos ^{2}v\right).
\end{eqnarray*}

By using these coefficients, the first fundamental form of $C$ is given by
\begin{eqnarray*}
ds^{2} &=&\{-\left( 1+k_{1}m_{1}p\left( u\right) \cos v-h^{\prime }\left(
u\right) \right) ^{2}+\left( k_{1}h\left( u\right) +k_{2}m_{2}p\left(
u\right) \sin v-m_{1}p^{\prime }\left( u\right) \cos v\right) ^{2} \\
&&+\left( k_{2}m_{1}p\left( u\right) \cos v+m_{2}p^{\prime }\left( u\right)
\sin v\right) ^{2}\}du^{2} \\
&& \\
&&+2\{p\left( u\right) \left( k_{1}m_{1}h\left( u\right) \sin
v+k_{2}m_{1}m_{2}p\left( u\right) +\left( -m_{1}^{2}+m_{2}^{2}\right)
p^{\prime }\left( u\right) \cos v\sin v\right) \}dudv \\
&& \\
&&+\{p^{2}\left( u\right) \left( m_{1}^{2}\sin ^{2}v+m_{2}^{2}\cos
^{2}v\right) \}dv^{2}.
\end{eqnarray*}

At the point $C(u,v)$, we have the following Lorentzian space-like angle:
\begin{eqnarray*}
\cos \psi &=&\dfrac{\{-\left( 1+k_{1}m_{1}p\left( u\right) \cos v-h^{\prime }\left(
u\right) \right) ^{2}+\left( k_{1}h\left( u\right) +k_{2}m_{2}p\left(
u\right) \sin v-m_{1}p^{\prime }\left( u\right) \cos v\right) ^{2}}{\sqrt{%
\Delta }} \\
&&\dfrac{+\left( k_{2}m_{1}p\left( u\right) \cos v+m_{2}p^{\prime }\left(
u\right) \sin v\right) ^{2}\}}{\sqrt{\Delta }}du \\
&& \\
&&+\dfrac{\{p\left( u\right) \left( k_{1}m_{1}h\left( u\right) \sin
v+k_{2}m_{1}m_{2}p\left( u\right) +\left( -m_{1}^{2}+m_{2}^{2}\right)
p^{\prime }\left( u\right) \cos v\sin v\right) \}}{\sqrt{\Delta }}dv,
\end{eqnarray*}
where
\begin{eqnarray*}
\Delta  &=&\{-\left( 1+k_{1}m_{1}p\left( u\right) \cos v-h^{\prime }\left( u\right)
\right) ^{2}+\left( k_{1}h\left( u\right) +k_{2}m_{2}p\left( u\right) \sin
v-m_{1}p^{\prime }\left( u\right) \cos v\right) ^{2} \\
&&+\left( k_{2}m_{1}p\left( u\right) \cos v+m_{2}p^{\prime }\left( u\right)
\sin v\right) ^{2}\}^{2}du \\
&& \\
&&+\{2p\left( u\right) \left( k_{1}m_{1}h\left( u\right) \sin
v+k_{2}m_{1}m_{2}p\left( u\right) +\left( -m_{1}^{2}+m_{2}^{2}\right)
p^{\prime }\left( u\right) \cos v\sin v\right) \\
&&\times\lbrack -\left( 1+k_{1}m_{1}p\left( u\right) \cos v-h^{\prime }\left(
u\right) \right) ^{2}+\left( k_{1}h\left( u\right) +k_{2}m_{2}p\left(
u\right) \sin v-m_{1}p^{\prime }\left( u\right) \cos v\right) ^{2} \\
&&+\left( k_{2}m_{1}p\left( u\right) \cos v+m_{2}p^{\prime }\left( u\right)
\sin v\right) ^{2}]\}dudv \\
&& \\
&&+\{p^{2}\left( u\right) \left( m_{1}^{2}\sin ^{2}v+m_{2}^{2}\cos
^{2}v\right) [-\left( 1+k_{1}m_{1}p\left( u\right) \cos v-h^{\prime }\left(
u\right) \right) ^{2} \\
&&+\left( k_{1}h\left( u\right) +k_{2}m_{2}p\left( u\right) \sin
v-m_{1}p^{\prime }\left( u\right) \cos v\right) ^{2}+\left(
k_{2}m_{1}p\left( u\right) \cos v+m_{2}p^{\prime }\left( u\right) \sin
v\right) ^{2}]\}dv^{2}.
\end{eqnarray*}
By using this equation, we find the following differential equation of the loxodrome:
\begin{eqnarray*}
&&\{-p^{2}\left( u\right) (k_{1}m_{1}h\left( u\right) \sin
v+k_{2}m_{1}m_{2}p(u)+\left( -m_{1}^{2}+m_{2}^{2}\right) p^{\prime }\left(
u\right) \cos v\sin v)^{2} \\
&&+\cos ^{2}\psi p^{2}\left( u\right) \left( m_{1}^{2}\sin
^{2}v+m_{2}^{2}\cos ^{2}v\right) [-(1+k_{1}m_{1}p\left( u\right) \cos
v-h^{\prime }\left( u\right) )^{2} \\
&&+\left( k_{1}h\left( u\right) +k_{2}m_{2}p(u)\sin v-m_{1}p^{\prime }\left(
u\right) \cos v\right) ^{2}+\left( k_{2}m_{1}p\left( u\right) \cos
v+m_{2}p^{\prime }\left( u\right) \sin v\right) ^{2}]\}\left( \dfrac{dv}{du}%
\right) ^{2} \\
&&-\{2\sin ^{2}\psi p\left( u\right) \left( k_{1}m_{1}h\left( u\right) \sin
v+k_{2}m_{1}m_{2}p\left( u\right) +\left( -m_{1}^{2}+m_{2}^{2}\right)
p^{\prime }\left( u\right) \cos v\sin v\right) \\
&&\times[-\left( 1+k_{1}m_{1}p\left( u\right) \cos v-h^{\prime }\left( u\right)
\right) ^{2}+\left( k_{1}h\left( u\right) +k_{2}m_{2}p\left( u\right) \sin
v-m_{1}p^{\prime }\left( u\right) \cos v\right) ^{2} \\
&&+\left( k_{2}m_{1}p\left( u\right) \cos v+m_{2}p^{\prime }\left( u\right)
\sin v\right) ^{2}]\}\dfrac{dv}{du} \\
&& \\
&=&\sin ^{2}\psi \{-\left( 1+k_{1}m_{1}p\left( u\right) \cos v-h^{\prime
}\left( u\right) \right) ^{2}+\left( k_{1}h\left( u\right)
+k_{2}m_{2}p\left( u\right) \sin v-m_{1}p^{\prime }\left( u\right) \cos
v\right) ^{2} \\
&&+\left( k_{2}m_{1}p\left( u\right) \cos v+m_{2}p^{\prime }\left( u\right)
\sin v\right) ^{2}\}^{2}.
\end{eqnarray*}

\textbf{4.} Let us consider the following space-like canal surface $C$ which is given by (3.41) in \cite{ucum}:

\begin{eqnarray}
C\left( u,v\right) &=&\alpha \left( u\right) +h\left( u\right) T\left(
u\right) +b(u,v)N\left( u\right) +\dfrac{t\left( u\right) }{2b(u,v)}B\left(
u\right)  \notag \\
&=&\alpha \left( u\right) +h\left( u\right) T\left( u\right) +bN\left(
u\right) +\dfrac{t\left( u\right) }{2b}B\left( u\right),
\end{eqnarray}
where $h\left( u\right) =r\left( u\right) r^{\prime }\left( u\right)$ and $%
t\left( u\right) =-r^{2}\left( u\right) \left( 1+r^{\prime 2}\left( u\right)
\right)$.

The coefficients of first fundamental form of the canal surface $C$ are
computed as
\begin{eqnarray*}
E &=&\left( 1-\frac{k_{1}t\left( u\right) }{2b}%
+h^{\prime }\left( u\right) \right) ^{2}\_\frac{\left(
b_{u}+k_{2}b+k_{1}h\left( u\right) \right) \left( b_{u}t\left( u\right)
+k_{2}bt\left( u\right) -bt^{\prime }\left( u\right) \right) }{b^{2}}, \\
&&\text{ \ \ \ \ \ \ \ \ \ \ \ \ } \\
F &=&-\frac{b_{v}\left( 2b_{u}t\left( u\right)
+k_{2}bt\left( u\right) +k_{1}h\left( u\right) t\left( u\right)
+k_{2}bt\left( u\right) -bt^{\prime }\left( u\right) \right) }{2b^{2}}, \\
&&\text{\ \ \ \ \ \ \ \ \ \ \ \ \ } \\
G &=&-\frac{b_{v}^{2}t\left( u\right) }{b^{2}},
\end{eqnarray*}
where $b_{u}$, $b_{v}$ refer to the derivative of the functions with respect to $u$, $v,$ respectively.

By using these coefficients, the first fundamental form of $C$\ is given by
\begin{eqnarray*}
ds^{2} &=&\{\left( 1-\frac{k_{1}t\left( u\right) }{2b}+h^{\prime }\left(
u\right) \right) ^{2}\_\frac{\left( b_{u}+k_{2}b+k_{1}h\left( u\right)
\right) \left( b_{u}t\left( u\right) +k_{2}bt\left( u\right) -bt^{\prime
}\left( u\right) \right) }{b^{2}}\}du^{2} \\
&& \\
&&-\{\frac{b_{v}\left( 2b_{u}t\left( u\right) +k_{2}bt\left( u\right)
+k_{1}h\left( u\right) t\left( u\right) +k_{2}bt\left( u\right) -bt^{\prime
}\left( u\right) \right) }{b^{2}}\}dudv-\frac{b_{v}^{2}t\left( u\right) }{b^{2}}dv^{2}.
\end{eqnarray*}
At the point $C(u,v)$, where the space-like loxodrome cuts the space-like meridians at a constant
Lorentzian space-like angle, we obtain

\begin{eqnarray*}
\cos \psi &=&\dfrac{\{\left( 1-\dfrac{k_{1}t\left( u\right) }{2b}+h^{\prime }\left(
u\right) \right) ^{2}\_\dfrac{\left( b_{u}+k_{2}b+k_{1}h\left( u\right)
\right) \left( b_{u}t\left( u\right) +k_{2}bt\left( u\right) -bt^{\prime
}\left( u\right) \right) }{b^{2}}\}}{\sqrt{\Delta }}du \\
&& \\
&&-\dfrac{\{\dfrac{b_{v}\left( 2b_{u}t\left( u\right) +k_{2}bt\left(
u\right) +k_{1}h\left( u\right) t\left( u\right) +k_{2}bt\left( u\right)
-bt^{\prime }\left( u\right) \right) }{2b^{2}}\}}{\sqrt{\Delta }}dv,
\end{eqnarray*}
where
\begin{eqnarray*}
\Delta &=&\{\left( 1-\dfrac{k_{1}t\left( u\right) }{2b}+h^{\prime }\left( u\right)
\right) ^{2}\_\dfrac{\left( b_{u}+k_{2}b+k_{1}h\left( u\right) \right)
\left( b_{u}t\left( u\right) +k_{2}bt\left( u\right) -bt^{\prime }\left(
u\right) \right) }{b^{2}}\}^{2}du \\
&& \\
&&-\{\frac{1}{b^{2}}b_{v}\left( 2b_{u}t\left( u\right) +k_{2}bt\left(
u\right) +k_{1}h\left( u\right) t\left( u\right) +k_{2}bt\left( u\right)
-bt^{\prime }\left( u\right) \right) \\
&&\times[\left( 1-\dfrac{k_{1}t\left( u\right) }{2b}+h^{\prime }\left( u\right)
\right) ^{2}-\dfrac{\left( b_{u}+k_{2}b+k_{1}h\left( u\right) \right) \left(
b_{u}t\left( u\right) +k_{2}bt\left( u\right) -bt^{\prime }\left( u\right)
\right) }{b^{2}}]\}dudv \\
&& \\
&&-\{\frac{b_{v}^{2}t\left( u\right) \left[ \left( 1-\dfrac{k_{1}t\left(
u\right) }{2b}+h^{\prime }\left( u\right) \right) ^{2}-\dfrac{\left(
b_{u}+k_{2}b+k_{1}h\left( u\right) \right) \left( b_{u}t\left( u\right)
+k_{2}bt\left( u\right) -bt^{\prime }\left( u\right) \right) }{b^{2}}\right]
}{b^{2}}\}dv^{2}. \\
&&
\end{eqnarray*}
From here, we get the following differential equation:
\begin{eqnarray*}
&&\{\frac{1}{4b^{4}}b_{v}^{2}[-\left( 2b_{u}t\left( u\right) +k_{2}bt\left(
u\right) +k_{1}h\left( u\right) t\left( u\right) +k_{2}bt\left( u\right)
-bt^{\prime }\left( u\right) \right) ^{2} \\
&&-4b^{2}\cos ^{2}\psi t\left( u\right) (\left( 1-\dfrac{k_{1}t\left(
u\right) }{2b}+h^{\prime }\left( u\right) \right) ^{2}-\dfrac{\left(
b_{u}+k_{2}b+k_{1}h\left( u\right) \right) \left( b_{u}t\left( u\right)
+k_{2}bt\left( u\right) -bt^{\prime }\left( u\right) \right) }{b^{2}}%
)]\}\left( \dfrac{dv}{du}\right) ^{2} \\
&& \\
&&-\{-\frac{1}{b^{2}}b_{v}\sin ^{2}\psi \left( 2b_{u}t\left( u\right)
+k_{2}bt\left( u\right) +k_{1}h\left( u\right) t\left( u\right)
+k_{2}bt\left( u\right) -bt^{\prime }\left( u\right) \right) \\
&&\times\lbrack \left( 1-\dfrac{k_{1}t\left( u\right) }{2b}+h^{\prime }\left(
u\right) \right) ^{2}-\dfrac{\left( b_{u}+k_{2}b+k_{1}h\left( u\right)
\right) \left( b_{u}t\left( u\right) +k_{2}bt\left( u\right) -bt^{\prime
}\left( u\right) \right) }{b^{2}}]\}\dfrac{dv}{du} \\
&& \\
&=&\sin ^{2}\psi \{\left( 1-\dfrac{k_{1}t\left( u\right) }{2b}+h^{\prime
}\left( u\right) \right) ^{2}\_\dfrac{\left( b_{u}+k_{2}b+k_{1}h\left(
u\right) \right) \left( b_{u}t\left( u\right) +k_{2}bt\left( u\right)
-bt^{\prime }\left( u\right) \right) }{b^{2}}\}^{2}.
\end{eqnarray*}

\textbf{5.} Let us consider the following space-like canal surface $C$ which is given by (3.62) in \cite{ucum}:
\begin{eqnarray}
C\left( u,v\right) &=&\alpha \left( u\right) -\dfrac{r^{2}\left( u\right)
+b^{2}(u,v)}{2h\left( u\right) }T\left( u\right) +b(u,v)N\left( u\right)
+h\left( u\right) B\left( u\right)  \notag \\
&=&\alpha \left( u\right) -\dfrac{r^{2}\left( u\right) +b^{2}}{2h\left(
u\right) }T\left( u\right) +bN\left( u\right) +h\left( u\right) B\left(
u\right),
\end{eqnarray}
where $h\left( u\right) =r\left( u\right) r^{\prime }\left( u\right).$

The coefficients of first fundamental form of the canal surface $C$ are given by
\begin{eqnarray*}
E &=&\left( b_{u}-\frac{\left( b^{2}+r^{2}\left(
u\right) \right) k_{1}}{2h\left( u\right) }-k_{2}h\left( u\right) \right)
^{2}+2\left( 1+k_{2}b-\frac{2h\left( u\right) \left( h\left( u\right)
+b_{u}b\right) -\left( b^{2}+r^{2}\left( u\right) \right) h^{\prime }\left(
u\right) }{2h^{2}\left( u\right) }\right), \\
&&\text{ \ \ \ \ \ \ \ \ \ \ \ \ } \\
F &=&b_{v}\left( b_{u}-\frac{k_{1}\left(
b^{2}+r^{2}\left( u\right) \right) }{2h\left( u\right) }-k_{2}h\left(
u\right) \right) +\frac{b_{v}b\left( bk_{1}-h^{\prime }\left( u\right)
\right) }{h\left( u\right)}, \\
&&\text{\ \ \ \ \ \ \ \ \ \ \ \ \ } \\
G &=&b_{v}^{2},
\end{eqnarray*}
where $b_{u}$ and $b_{v}$ refer to the derivative of the functions with respect to $u$, $v,$ respectively.

By using these coefficients, the first fundamental form of $C$\ is given by
\begin{eqnarray*}
ds^{2} &=&\{\left( b_{u}-\frac{\left( b^{2}+r^{2}\left( u\right) \right)
k_{1}}{2h\left( u\right) }-k_{2}h\left( u\right) \right) ^{2}+2\left(
1+k_{2}b-\frac{2h\left( u\right) \left( h\left( u\right) +b_{u}b\right)
-\left( b^{2}+r^{2}\left( u\right) \right) h^{\prime }\left( u\right) }{%
2h^{2}\left( u\right) }\right) \}du^{2} \\
&& \\
&&+2\{b_{v}\left( b_{u}-\frac{k_{1}\left( b^{2}+r^{2}\left( u\right) \right)
}{2h\left( u\right) }-k_{2}h\left( u\right) \right) +\frac{b_{v}b\left(
bk_{1}-h^{\prime }\left( u\right) \right) }{h\left( u\right) }\}dudv+b_{v}^{2}dv^{2}.
\end{eqnarray*}

As it was mentioned earlier, at the point $C(u,v)$, where the space-like loxodrome cuts the space-like meridians at a
constant Lorentzian space-like angle, we have
\begin{eqnarray*}
\cos \psi &=&\dfrac{\{\left( b_{u}-\dfrac{\left( b^{2}+r^{2}\left( u\right) \right)
k_{1}}{2h\left( u\right) }-k_{2}h\left( u\right) \right) ^{2}+2\left(
1+k_{2}b-\dfrac{2h\left( u\right) \left( h\left( u\right) +b_{u}b\right)
-\left( b^{2}+r^{2}\left( u\right) \right) h^{\prime }\left( u\right) }{%
2h^{2}\left( u\right) }\right) \}}{\sqrt{\Delta }}du \\
&& \\
&&+\dfrac{\{b_{v}\left( b_{u}-\dfrac{k_{1}\left( b^{2}+r^{2}\left( u\right)
\right) }{2h\left( u\right) }-k_{2}h\left( u\right) \right) +\dfrac{%
b_{v}b\left( bk_{1}-h^{\prime }\left( u\right) \right) }{h\left( u\right) }\}%
}{\sqrt{\Delta }}dv,
\end{eqnarray*}
where
\begin{eqnarray*}
\Delta  &=&\{\left( b_{u}-\frac{\left( b^{2}+r^{2}\left( u\right) \right)
k_{1}}{2h\left( u\right) }-k_{2}h\left( u\right) \right) ^{2}+2\left(
1+k_{2}b-\frac{2h\left( u\right) \left( h\left( u\right) +b_{u}b\right)
-\left( b^{2}+r^{2}\left( u\right) \right) h^{\prime }\left( u\right) }{%
2h^{2}\left( u\right) }\right) \}^{2}du \\
&& \\
&&+\{2\left[ b_{v}\left( b_{u}-\frac{\left( b^{2}+r^{2}\left( u\right)
\right) k_{1}}{2h\left( u\right) }-k_{2}h\left( u\right) \right) +\frac{%
bb_{v}\left( bk_{1}-h^{\prime }\left( u\right) \right) }{h\left( u\right) }%
\right]\times\left[ \left( b_{u}-\frac{\left( b^{2}+r^{2}\left( u\right) \right)
k_{1}}{2h\left( u\right) }-k_{2}h\left( u\right) \right) ^{2}\right.  \\
&&\left. +2\left( 1+k_{2}b-\frac{2h\left( u\right) \left( h\left( u\right)
+b_{u}b\right) -\left( b^{2}+r^{2}\left( u\right) \right) h^{\prime }\left(
u\right) }{2h^{2}\left( u\right) }\right) \right] \}dudv \\
&& \\
&&+\{b_{v}^{2}\left[ \left( b_{u}-\frac{\left( b^{2}+r^{2}\left( u\right)
\right) k_{1}}{2h\left( u\right) }-k_{2}h\left( u\right) \right) ^{2}\right.
\\
&&\left. +2\left( 1+k_{2}b-\frac{2h\left( u\right) \left( h\left( u\right)
+b_{u}b\right) -\left( b^{2}+r^{2}\left( u\right) \right) h^{\prime }\left(
u\right) }{2h^{2}\left( u\right) }\right) \right] \}dv^{2}
\end{eqnarray*}%
Thus the following differential equation is obtained:
\begin{eqnarray*}
&&\{-\left[ b_{v}\left( b_{u}-\dfrac{\left( b^{2}+r^{2}\left( u\right)
\right) k_{1}}{2h\left( u\right) }-k_{2}h\left( u\right) \right) +\frac{%
bb_{v}\left( bk_{1}-h^{\prime }\left( u\right) \right) }{h\left( u\right) }%
\right] ^{2} \\
&&+\cos ^{2}\psi b_{v}^{2}\left[ \left( b_{u}-\dfrac{\left(
b^{2}+r^{2}\left( u\right) \right) k_{1}}{2h\left( u\right) }-k_{2}h\left(
u\right) \right) ^{2}\right.  \\
&&\left. +2\left( 1+k_{2}b-\dfrac{2h\left( u\right) \left( h\left( u\right)
+b_{u}b\right) -\left( b^{2}+r^{2}\left( u\right) \right) h^{\prime }\left(
u\right) }{2h^{2}\left( u\right) }\right) \right] \}\left( \dfrac{dv}{du}%
\right) ^{2} \\
&& \\
&&-\{2\sin ^{2}\psi \left[ b_{v}\left( b_{u}-\dfrac{\left( b^{2}+r^{2}\left(
u\right) \right) k_{1}}{2h\left( u\right) }-k_{2}h\left( u\right) \right) +%
\frac{bb_{v}\left( bk_{1}-h^{\prime }\left( u\right) \right) }{h\left(
u\right) }\right]  \\
&&\times\left[ \left( b_{u}-\dfrac{\left( b^{2}+r^{2}\left( u\right) \right) k_{1}%
}{2h\left( u\right) }-k_{2}h\left( u\right) \right) ^{2}+2\left( 1+k_{2}b-%
\dfrac{2h\left( u\right) \left( h\left( u\right) +b_{u}b\right) -\left(
b^{2}+r^{2}\left( u\right) \right) h^{\prime }\left( u\right) }{2h^{2}\left(
u\right) }\right) \right] \}\dfrac{dv}{du} \\
&& \\
&=&\sin ^{2}\psi \{\left( b_{u}-\dfrac{\left( b^{2}+r^{2}\left( u\right)
\right) k_{1}}{2h\left( u\right) }-k_{2}h\left( u\right) \right)
^{2}+2\left( 1+k_{2}b-\dfrac{2h\left( u\right) \left( h\left( u\right)
+b_{u}b\right) -\left( b^{2}+r^{2}\left( u\right) \right) h^{\prime }\left(
u\right) }{2h^{2}\left( u\right) }\right) \}^{2}.
\end{eqnarray*}

\section{Differential equations of the space-like loxodromes on the
time-like canal surfaces having space-like meridians}

In this section, we investigate the differential equations of space-like
loxodromes on the time-like canal surfaces having space-like meridians (that is, $EG-F^{2}<0$ and $E>0$ for all $(u,v)$).

\begin{definition}
\label{d6} If $u$ and $v$ are space-like vectors in ${\mathbb{E}}_{1}^{3}$
which span a time-like plane. Then
\begin{equation}
|\langle u,v\rangle| =\Vert u\Vert \Vert v\Vert \cosh \eta ,  \notag
\end{equation}%
where $\eta \in {\mathbb{R}^{+}}$ is the Lorentzian
time-like angle between $u$ and $v$ \cite{ratcliffe}.
\end{definition}
\begin{definition}
\label{d7} A space-like curve on a time-like canal surface in ${\mathbb{E}}%
^3_1$ having space-like meridians is called as a loxodrome if the curve cuts
all meridians at a constant Lorentzian time-like angle.
\end{definition}

\textbf{1. }Let us consider the following time-like canal surface $C$ which is given by (3.58) in \cite{ucum}:
\begin{eqnarray}
C\left( u,v\right) &=&\alpha \left( u\right) -h\left( u\right) T\left(
u\right) +b(u,v)N\left( u\right) +\dfrac{p\left( u\right) }{2b(u,v)}B\left(
u\right)  \notag \\
&=&\alpha \left( u\right) -h\left( u\right) T\left( u\right) +bN\left(
u\right) +\dfrac{p\left( u\right) }{2b}B\left( u\right),
\end{eqnarray}
where
$h\left( u\right) =r\left( u\right) r^{\prime }\left(u\right)$ and $p\left( u\right) =r^{2}\left( u\right) \left( 1-r^{\prime2}\left( u\right) \right)$.

The coefficients of first fundamental form of the canal surface $C$ are
\begin{eqnarray*}
E &=&\left( -1+\frac{k_{1}p\left( u\right) }{2b}%
+h^{\prime }\left( u\right) \right) ^{2}-\frac{\left(
b_{u}+k_{2}b-k_{1}h\left( u\right) \right) [p\left( u\right) \left(
b_{u}+k_{2}b\right) -p^{\prime }\left( u\right) b]}{b^{2}}, \\
&&\text{ \ \ \ \ \ \ \ \ \ \ \ \ } \\
F &=&-\frac{b_{v}[p\left( u\right) \left(
2b_{u}+k_{2}b-k_{1}h\left( u\right) +k_{2}b\right) -p^{\prime }\left(
u\right) b]}{2b^{2}}, \\
&&\text{\ \ \ \ \ \ \ \ \ \ \ \ \ } \\
G &=&-\dfrac{p\left( u\right) b_{v}^{2}}{b^{2}}.
\end{eqnarray*}
By using these equations, the first fundamental form of $C$ is
\begin{eqnarray*}
ds^{2} &=&\{\left( -1+\frac{k_{1}p\left( u\right) }{2b}+h^{\prime }\left(
u\right) \right) ^{2}-\frac{\left( b_{u}+k_{2}b-k_{1}h\left( u\right)
\right) [p\left( u\right) \left( b_{u}+k_{2}b\right) -p^{\prime }\left(
u\right) b]}{b^{2}}\}du^{2} \\
&& \\
&&+\{-\frac{b_{v}[p\left( u\right) \left( 2b_{u}+k_{2}b-k_{1}h\left(
u\right) +k_{2}b\right) -p^{\prime }\left( u\right) b]}{b^{2}}%
\}dudv-\dfrac{p\left( u\right) b_{v}^{2}}{b^{2}}dv^{2}.
\end{eqnarray*}
Lorentzian time-like angle between the space-like loxodrome and the space-like meridian is defined by the angle between their tangent vectors
at the point $C(u,v)$ that can be given by
\begin{eqnarray*}
\epsilon \cosh \eta &=&\dfrac{Edu+Fdv}{\sqrt{\Gamma }} \\
&& \\
&=&\dfrac{\{\left( -1+\dfrac{k_{1}p\left( u\right) }{2b}+h^{\prime }\left(
u\right) \right) ^{2}-\dfrac{\left( b_{u}+k_{2}b-k_{1}h\left( u\right)
\right) [p\left( u\right) \left( b_{u}+k_{2}b\right) -p^{\prime }\left(
u\right) b]}{b^{2}}\}}{\sqrt{\Gamma }}du \\
&& \\
&&+\dfrac{\{-\dfrac{b_{v}[p\left( u\right) \left( 2b_{u}+k_{2}b-k_{1}h\left(
u\right) +k_{2}b\right) -p^{\prime }\left( u\right) b]}{2b^{2}}\}}{\sqrt{%
\Gamma }}dv,
\end{eqnarray*}
where $\epsilon =\pm 1$ and
\begin{eqnarray*}
\Gamma &=&E^{2}du^{2}+2EFdudv+EGdv^{2} \\
&& \\
&=&\{\left( -1+\dfrac{k_{1}p\left( u\right) }{2b}+h^{\prime }\left( u\right)
\right) ^{2}-\dfrac{\left( b_{u}+k_{2}b-k_{1}h\left( u\right) \right)
[p\left( u\right) \left( b_{u}+k_{2}b\right) -p^{\prime }\left( u\right) b]}{%
b^{2}}\}^{2}du^{2} \\
&& \\
&&+\{-\frac{1}{b^{2}}b_{v}[p\left( u\right) \left(
2b_{u}+k_{2}b-k_{1}h\left( u\right) +k_{2}b\right) -p^{\prime }\left(
u\right) b] \\
&&\times[\left( -1+\dfrac{k_{1}p\left( u\right) }{2b}+h^{\prime }\left( u\right)
\right) ^{2}-\dfrac{\left( b_{u}+k_{2}b-k_{1}h\left( u\right) \right)
[p\left( u\right) \left( b_{u}+k_{2}b\right) -p^{\prime }\left( u\right) b]}{%
b^{2}}]\}dudv \\
&& \\
&&+\{-\frac{b_{v}^{2}p\left( u\right) \left[ \left( -1+\dfrac{k_{1}p\left(
u\right) }{2b}+h^{\prime }\left( u\right) \right) ^{2}-\dfrac{\left(
b_{u}+k_{2}b-k_{1}h\left( u\right) \right) [p\left( u\right) \left(
b_{u}+k_{2}b\right) -p^{\prime }\left( u\right) b]}{b^{2}}\right] }{b^{2}}%
\}dv^{2}.
\end{eqnarray*}
From this equation, we obtain the following differential equation of space-like loxodromes on the time-like canal surfaces having space-like meridians:
\begin{equation*}
\left( -\cosh ^{2}\eta EG+F^{2}\right) \left( \dfrac{dv}{du}\right)
^{2}-2\sinh ^{2}\eta EF\dfrac{dv}{du}=\sinh ^{2}\eta E^{2},
\end{equation*}%
that is
\begin{eqnarray*}
&&\{\frac{1}{4b^{4}}b_{v}^{2}\text{ }[(p\left( u\right) \left(
2b_{u}+k_{2}b-k_{1}h\left( u\right) +k_{2}b\right) -p^{\prime }\left(
u\right) b)^{2}+4\cosh ^{2}\eta b^{2}p\left( u\right) \\
&&\times(\left( -1+\dfrac{k_{1}p\left( u\right) }{2b}+h^{\prime }\left( u\right)
\right) ^{2}-\dfrac{\left( b_{u}+k_{2}b-k_{1}h\left( u\right) \right)
[p\left( u\right) \left( b_{u}+k_{2}b\right) -p^{\prime }\left( u\right) b]}{%
b^{2}})]\}\left( \dfrac{dv}{du}\right) ^{2} \\
&& \\
&&+\{\frac{1}{b^{2}}\sinh ^{2}\eta b_{v}[p\left( u\right) \left(
2b_{u}+k_{2}b-k_{1}h\left( u\right) +k_{2}b\right) -p^{\prime }\left(
u\right) b] \\
&&\times[ \left( -1+\dfrac{k_{1}p\left( u\right) }{2b}+h^{\prime }\left(
u\right) \right) ^{2}-\dfrac{\left( b_{u}+k_{2}b-k_{1}h\left( u\right)
\right) [p\left( u\right) \left( b_{u}+k_{2}b\right) -p^{\prime }\left(
u\right) b]}{b^{2}}]\}\dfrac{dv}{du} \\
&& \\
&=&\sinh ^{2}\eta \{\left( -1+\dfrac{k_{1}p\left( u\right) }{2b}+h^{\prime
}\left( u\right) \right) ^{2}-\dfrac{\left( b_{u}+k_{2}b-k_{1}h\left(
u\right) \right) [p\left( u\right) \left( b_{u}+k_{2}b\right) -p^{\prime
}\left( u\right) b]}{b^{2}}\}^{2}.
\end{eqnarray*}

\textbf{2.} Let us consider the following time-like canal surface $C$ which is given by (3.65) in \cite{ucum}:
\begin{eqnarray}
C\left( u,v\right) &=&\alpha \left( u\right) +\dfrac{b^{2}(u,v)-r^{2}\left(
u\right) }{2h\left( u\right) }T\left( u\right) +b(u,v)N\left( u\right)
-h\left( u\right) B\left( u\right)  \notag \\
&=&\alpha \left( u\right) +\dfrac{b^{2}-r^{2}\left( u\right) }{2h\left(
u\right) }T\left( u\right) +bN\left( u\right) -h\left( u\right) B\left(
u\right),
\end{eqnarray}
where $h\left( u\right) =r\left( u\right) r^{\prime }\left(
u\right)$.

The coefficients of first fundamental form of the canal surface $C$ are
computed as
\begin{eqnarray*}
E &=&\left( b_{u}+\frac{k_{1}\left(
b^{2}-r^{2}\left( u\right) \right) }{2h\left( u\right) }+k_{2}h\left(
u\right) \right) ^{2} \\
&&+\frac{\left( k_{1}b+h^{\prime }\left( u\right) \right) [2h^{2}\left(
u\right) -2b_{u}bh\left( u\right) -2h^{2}\left( u\right) \left(
1+k_{2}b\right) +\left( b^{2}-r^{2}\left( u\right) \right) h^{\prime }\left(
u\right) ]}{h^{2}\left( u\right), } \\
&&\text{\ \ \ \ \ \ \ \ \ } \\
F &=&b_{v}\left( b_{u}+\frac{k_{1}\left(
b^{2}-r^{2}\left( u\right) \right) }{2h\left( u\right) }+k_{2}h\left(
u\right) \right) +\frac{b_{v}b\left( -k_{1}b+h^{\prime }\left( u\right)
\right) }{h\left( u\right)}, \\
&&\text{\ \ \ \ \ \ \ \ \ \ \ \ \ } \\
G &=&b_{v}^{2},
\end{eqnarray*}
where $b_{u}$, $b_{v}$ refer to the derivative of the functions with respect to $u$, $v,$ respectively.

By using these equations, the first fundamental form of $C$\ is given by
\begin{eqnarray*}
ds^{2} &=&\{\left( b_{u}+\frac{k_{1}\left( b^{2}-r^{2}\left( u\right)
\right) }{2h\left( u\right) }+k_{2}h\left( u\right) \right) ^{2} \\
&&+\frac{\left( k_{1}b+h^{\prime }\left( u\right) \right) [2h^{2}\left(
u\right) -2b_{u}bh\left( u\right) -2h^{2}\left( u\right) \left(
1+k_{2}b\right) +\left( b^{2}-r^{2}\left( u\right) \right) h^{\prime }\left(
u\right) ]}{h^{2}\left( u\right) }\}du^{2} \\
&& \\
&&+2\{b_{v}\left( b_{u}+\frac{k_{1}\left( b^{2}-r^{2}\left( u\right) \right)
}{2h\left( u\right) }+k_{2}h\left( u\right) \right) +\frac{b_{v}b\left(
-k_{1}b+h^{\prime }\left( u\right) \right) }{h\left( u\right) }%
\}dudv+\left\{ b_{v}^{2}\right\} dv^{2}.
\end{eqnarray*}
At the point $C(u,v)$, where the space-like loxodrome cuts the space-like meridians at a constant
Lorentzian time-like angle, we have
\begin{eqnarray*}
\epsilon \cosh \eta &=&\dfrac{\{\left( b_{u}+\dfrac{k_{1}\left( b^{2}-r^{2}\left( u\right)
\right) }{2h\left( u\right) }+k_{2}h\left( u\right) \right) ^{2}}{\sqrt{%
\Gamma }} \\
&&+\frac{\dfrac{\left( k_{1}b+h^{\prime }\left( u\right) \right)
[2h^{2}\left( u\right) -2b_{u}bh\left( u\right) -2h^{2}\left( u\right)
\left( 1+k_{2}b\right) +\left( b^{2}-r^{2}\left( u\right) \right) h^{\prime
}\left( u\right) ]}{h^{2}\left( u\right) }\}}{\sqrt{\Gamma }}du \\
&& \\
&&+\dfrac{\{b_{v}\left( b_{u}+\dfrac{k_{1}\left( b^{2}-r^{2}\left( u\right)
\right) }{2h\left( u\right) }+k_{2}h\left( u\right) \right) +\dfrac{%
b_{v}b\left( -k_{1}b+h^{\prime }\left( u\right) \right) }{h\left( u\right) }%
\}}{\sqrt{\Gamma }}dv,
\end{eqnarray*}
where $\epsilon =\pm 1$ and
\begin{eqnarray*}
\Gamma &=&\{\left( b_{u}+\dfrac{k_{1}\left( b^{2}-r^{2}\left( u\right) \right) }{%
2h\left( u\right) }+k_{2}h\left( u\right) \right) ^{2} \\
&&+\dfrac{\left( k_{1}b+h^{\prime }\left( u\right) \right) [2h^{2}\left(
u\right) -2b_{u}bh\left( u\right) -2h^{2}\left( u\right) \left(
1+k_{2}b\right) +\left( b^{2}-r^{2}\left( u\right) \right) h^{\prime }\left(
u\right) ]}{h^{2}\left( u\right) }\}^{2}du \\
&& \\
&&+\{2\left[ b_{v}\left( b_{u}+\dfrac{k_{1}\left( b^{2}-r^{2}\left( u\right)
\right) }{2h\left( u\right) }+k_{2}h\left( u\right) \right) +\frac{%
bb_{v}\left( -k_{1}b+h^{\prime }\left( u\right) \right) }{h\left( u\right) }%
\right] [\left( b_{u}+\dfrac{k_{1}\left( b^{2}-r^{2}\left( u\right) \right)
}{2h\left( u\right) }+k_{2}h\left( u\right) \right) ^{2} \\
&&+\dfrac{\left( k_{1}b+h^{\prime }\left( u\right) \right) [2h^{2}\left(
u\right) -2b_{u}bh\left( u\right) -2h^{2}\left( u\right) \left(
1+k_{2}b\right) +\left( b^{2}-r^{2}\left( u\right) \right) h^{\prime }\left(
u\right) ]}{h^{2}\left( u\right) }]\}dudv \\
&& \\
&&+\{b_{v}^{2}[\left( b_{u}+\dfrac{k_{1}\left( b^{2}-r^{2}\left( u\right)
\right) }{2h\left( u\right) }+k_{2}h\left( u\right) \right) ^{2} \\
&&+\dfrac{\left( k_{1}b+h^{\prime }\left( u\right) \right) [2h^{2}\left(
u\right) -2b_{u}bh\left( u\right) -2h^{2}\left( u\right) \left(
1+k_{2}b\right) +\left( b^{2}-r^{2}\left( u\right) \right) h^{\prime }\left(
u\right) ]}{h^{2}\left( u\right) }]\}dv^{2}.
\end{eqnarray*}
From this equation, we obtain the following differential equation:
\begin{eqnarray*}
&&\{\left[ b_{v}\left( b_{u}+\dfrac{k_{1}\left( b^{2}-r^{2}\left( u\right)
\right) }{2h\left( u\right) }+k_{2}h\left( u\right) \right) +\frac{%
bb_{v}\left( -k_{1}b+h^{\prime }\left( u\right) \right) }{h\left( u\right) }%
\right] ^{2} \\
&&-\cosh ^{2}\eta b_{v}^{2}\left[ \left( b_{u}+\dfrac{k_{1}\left(
b^{2}-r^{2}\left( u\right) \right) }{2h\left( u\right) }+k_{2}h\left(
u\right) \right) ^{2}\right.  \\
&&\left. +\dfrac{\left( k_{1}b+h^{\prime }\left( u\right) \right)
[2h^{2}\left( u\right) -2b_{u}bh\left( u\right) -2h^{2}\left( u\right)
\left( 1+k_{2}b\right) +\left( b^{2}-r^{2}\left( u\right) \right) h^{\prime
}\left( u\right) ]}{h^{2}\left( u\right) }\right] \}\left( \dfrac{dv}{du}%
\right) ^{2} \\
&& \\
&&-\{2\sinh ^{2}\eta \left[ b_{v}\left( b_{u}+\dfrac{k_{1}\left(
b^{2}-r^{2}\left( u\right) \right) }{2h\left( u\right) }+k_{2}h\left(
u\right) \right) +\frac{bb_{v}\left( -k_{1}b+h^{\prime }\left( u\right)
\right) }{h\left( u\right) }\right]  \\
&&\times\left[ \left( b_{u}+\dfrac{k_{1}\left( b^{2}-r^{2}\left( u\right) \right)
}{2h\left( u\right) }+k_{2}h\left( u\right) \right) ^{2}\right.  \\
&&\left. +\dfrac{\left( k_{1}b+h^{\prime }\left( u\right) \right)
[2h^{2}\left( u\right) -2b_{u}bh\left( u\right) -2h^{2}\left( u\right)
\left( 1+k_{2}b\right) +\left( b^{2}-r^{2}\left( u\right) \right) h^{\prime
}\left( u\right) ]}{h^{2}\left( u\right) }\right] \}\dfrac{dv}{du} \\
&& \\
&=&\sinh ^{2}\eta \{\left( b_{u}+\dfrac{k_{1}\left( b^{2}-r^{2}\left(
u\right) \right) }{2h\left( u\right) }+k_{2}h\left( u\right) \right) ^{2} \\
&&+\dfrac{\left( k_{1}b+h^{\prime }\left( u\right) \right) [2h^{2}\left(
u\right) -2b_{u}bh\left( u\right) -2h^{2}\left( u\right) \left(
1+k_{2}b\right) +\left( b^{2}-r^{2}\left( u\right) \right) h^{\prime }\left(
u\right) ]}{h^{2}\left( u\right) }\}^{2}.
\end{eqnarray*}

\section{Differential equations of the space-like loxodromes on the
time-like canal surfaces having time-like meridians}
In this section, we obtain the differential equations of space-like loxodromes on the time-like canal surfaces having time-like meridians (that is, $EG-F^{2}<0$ and $E<0$ for all $(u,v)$).
\begin{definition}
\label{d10}
If $u$ is a space-like vector and $v$ is a time-like vector in ${\mathbb{E}}^3_1$. Then
\begin{equation}
|\langle u,v\rangle|=\|u\|\|v\|\sinh\varphi,  \notag
\end{equation}
where $\varphi\in{\mathbb{R}^{+}}\cup \{0\}$ is the Lorentzian time-like
angle between $u$ and $v$ \cite{ratcliffe}.
\end{definition}
\begin{definition}
\label{d11}
A space-like curve on a time-like canal surface in ${\mathbb{E}}^3_1$ having
time-like meridians is called as a loxodrome if the curve cuts all
meridians at a constant Lorentzian time-like angle.
\end{definition}

\textbf{1. }Let us consider the following time-like canal surface $C$ which
is given by (3.58) in \cite{ucum}:
\begin{eqnarray}
C\left( u,v\right) &=&\alpha \left( u\right) -h\left( u\right) T\left(
u\right) +b(u,v)N\left( u\right) +\dfrac{p\left( u\right) }{2b(u,v)}B\left(
u\right)  \notag \\
&=&\alpha \left( u\right) -h\left( u\right) T\left( u\right) +bN\left(
u\right) +\dfrac{p\left( u\right) }{2b}B\left( u\right),
\end{eqnarray}
where
$h\left( u\right) =r\left( u\right) r^{\prime }\left(u\right) $ and $p\left( u\right) =r^{2}\left( u\right) \left( 1-r^{\prime2}\left( u\right) \right)$.

At the point $C(u,v)$, where the space-like loxodrome cuts the time-like meridians at a constant Lorentzian time-like angle, we have
\begin{eqnarray*}
\epsilon \sinh \varphi &=&\dfrac{Edu+Fdv}{\sqrt{\Omega }} \\
&& \\
&=&\dfrac{\{\left( -1+\dfrac{k_{1}p\left( u\right) }{2b}+h^{\prime }\left(
u\right) \right) ^{2}-\dfrac{\left( b_{u}+k_{2}b-k_{1}h\left( u\right)
\right) [p\left( u\right) \left( b_{u}+k_{2}b\right) -p^{\prime }\left(
u\right) b]}{b^{2}}\}}{\sqrt{\Omega }}du \\
&& \\
&&+\dfrac{\{-\dfrac{b_{v}[p\left( u\right) \left( 2b_{u}+k_{2}b-k_{1}h\left(
u\right) +k_{2}b\right) -p^{\prime }\left( u\right) b]}{2b^{2}}\}}{\sqrt{%
\Omega }}dv,
\end{eqnarray*}
where $\epsilon =\pm 1$ and
\begin{eqnarray*}
\Omega &=&-E^{2}du^{2}-2EFdudv-EGdv^{2} \\
&& \\
&=&-\{\left( -1+\dfrac{k_{1}p\left( u\right) }{2b}+h^{\prime }\left(
u\right) \right) ^{2}-\dfrac{\left( b_{u}+k_{2}b-k_{1}h\left( u\right)
\right) [p\left( u\right) \left( b_{u}+k_{2}b\right) -p^{\prime }\left(
u\right) b]}{b^{2}}\}^{2}du^{2} \\
&& \\
&&-\{-\frac{1}{b^{2}}b_{v}[p\left( u\right) \left(
2b_{u}+k_{2}b-k_{1}h\left( u\right) +k_{2}b\right) -p^{\prime }\left(
u\right) b] \\
&&\times[\left( -1+\dfrac{k_{1}p\left( u\right) }{2b}+h^{\prime }\left( u\right)
\right) ^{2}-\dfrac{\left( b_{u}+k_{2}b-k_{1}h\left( u\right) \right)
[p\left( u\right) \left( b_{u}+k_{2}b\right) -p^{\prime }\left( u\right) b]}{%
b^{2}}]\}dudv \\
&& \\
&&+\{\frac{b_{v}^{2}p\left( u\right) \left[ \left( -1+\dfrac{k_{1}p\left(
u\right) }{2b}+h^{\prime }\left( u\right) \right) ^{2}-\dfrac{\left(
b_{u}+k_{2}b-k_{1}h\left( u\right) \right) [p\left( u\right) \left(
b_{u}+k_{2}b\right) -p^{\prime }\left( u\right) b]}{b^{2}}\right] }{b^{2}}%
\}dv^{2}.
\end{eqnarray*}
From this equation, we obtain the following differential equation of space-like loxodromes on the time-like canal surfaces having time-like meridians:
\begin{equation*}
\left( \sinh^{2}\varphi EG+F^{2}\right) \left( \dfrac{dv}{du}\right)
^{2}+2\cosh ^{2}\varphi EF\dfrac{dv}{du}=-\cosh ^{2}\varphi E^{2},
\end{equation*}%
that is
\begin{eqnarray*}
&&\{\frac{1}{4b^{4}}b_{v}^{2}\text{ }[(p\left( u\right) \left(
2b_{u}+k_{2}b-k_{1}h\left( u\right) +k_{2}b\right) -p^{\prime }\left(
u\right) b)^{2}-4\sinh^{2}\varphi b^{2}p\left( u\right) \\
&&\times(\left( -1+\dfrac{k_{1}p\left( u\right) }{2b}+h^{\prime }\left( u\right)
\right) ^{2}-\dfrac{\left( b_{u}+k_{2}b-k_{1}h\left( u\right) \right)
[p\left( u\right) \left( b_{u}+k_{2}b\right) -p^{\prime }\left( u\right) b]}{%
b^{2}})]\}\left( \dfrac{dv}{du}\right) ^{2} \\
&& \\
&&+\{-\frac{1}{b^{2}}\cosh ^{2}\varphi b_{v}[p\left( u\right) \left(
2b_{u}+k_{2}b-k_{1}h\left( u\right) +k_{2}b\right) -p^{\prime }\left(
u\right) b] \\
&&\times[\left( -1+\dfrac{k_{1}p\left( u\right) }{2b}+h^{\prime }\left( u\right)
\right) ^{2}-\dfrac{\left( b_{u}+k_{2}b-k_{1}h\left( u\right) \right)
[p\left( u\right) \left( b_{u}+k_{2}b\right) -p^{\prime }\left( u\right) b]}{%
b^{2}}]\}\dfrac{dv}{du} \\
&& \\
&=&-\cosh ^{2}\varphi \{\left( -1+\dfrac{k_{1}p\left( u\right) }{2b}%
+h^{\prime }\left( u\right) \right) ^{2}-\dfrac{\left(
b_{u}+k_{2}b-k_{1}h\left( u\right) \right) [p\left( u\right) \left(
b_{u}+k_{2}b\right) -p^{\prime }\left( u\right) b]}{b^{2}}\}^{2}.
\end{eqnarray*}

\textbf{2.} Let us consider the following time-like canal surface $C$ which is given by (3.65) in \cite{ucum}:
\begin{eqnarray}
C\left( u,v\right) &=&\alpha \left( u\right) +\dfrac{b^{2}(u,v)-r^{2}\left(
u\right) }{2h\left( u\right) }T\left( u\right) +b(u,v)N\left( u\right)
-h\left( u\right) B\left( u\right)  \notag \\
&=&\alpha \left( u\right) +\dfrac{b^{2}-r^{2}\left( u\right) }{2h\left(
u\right) }T\left( u\right) +bN\left( u\right) -h\left( u\right) B\left(
u\right).
\end{eqnarray}
As it was mentioned earlier, at the point $C(u,v)$, where the space-like loxodrome cuts the time-like meridians at a constant Lorentzian time-like angle, we have
\begin{eqnarray*}
\epsilon \sinh \varphi &=&\dfrac{\{\left( b_{u}+\dfrac{k_{1}\left( b^{2}-r^{2}\left( u\right)
\right) }{2h\left( u\right) }+k_{2}h\left( u\right) \right) ^{2}}{\sqrt{%
\Omega }} \\
&&+\frac{\dfrac{\left( k_{1}b+h^{\prime }\left( u\right) \right)
[2h^{2}\left( u\right) -2b_{u}bh\left( u\right) -2h^{2}\left( u\right)
\left( 1+k_{2}b\right) +\left( b^{2}-r^{2}\left( u\right) \right) h^{\prime
}\left( u\right) ]}{h^{2}\left( u\right) }\}}{\sqrt{\Omega }}du \\
&& \\
&&+\dfrac{\{b_{v}\left( b_{u}+\dfrac{k_{1}\left( b^{2}-r^{2}\left( u\right)
\right) }{2h\left( u\right) }+k_{2}h\left( u\right) \right) +\dfrac{%
b_{v}b\left( -k_{1}b+h^{\prime }\left( u\right) \right) }{h\left( u\right) }%
\}}{\sqrt{\Omega }}dv,
\end{eqnarray*}
where $\epsilon =\pm 1$ and
\begin{eqnarray*}
\Omega &=&-\{\left( b_{u}+\dfrac{k_{1}\left( b^{2}-r^{2}\left( u\right) \right) }{%
2h\left( u\right) }+k_{2}h\left( u\right) \right) ^{2} \\
&&+\dfrac{\left( k_{1}b+h^{\prime }\left( u\right) \right) [2h^{2}\left(
u\right) -2b_{u}bh\left( u\right) -2h^{2}\left( u\right) \left(
1+k_{2}b\right) +\left( b^{2}-r^{2}\left( u\right) \right) h^{\prime }\left(
u\right) ]}{h^{2}\left( u\right) }\}^{2}du \\
&& \\
&&-\{2\left[ b_{v}\left( b_{u}+\dfrac{k_{1}\left( b^{2}-r^{2}\left( u\right)
\right) }{2h\left( u\right) }+k_{2}h\left( u\right) \right) +\frac{%
bb_{v}\left( -k_{1}b+h^{\prime }\left( u\right) \right) }{h\left( u\right) }%
\right] [\left( b_{u}+\dfrac{k_{1}\left( b^{2}-r^{2}\left( u\right) \right)
}{2h\left( u\right) }+k_{2}h\left( u\right) \right) ^{2} \\
&&+\dfrac{\left( k_{1}b+h^{\prime }\left( u\right) \right) [2h^{2}\left(
u\right) -2b_{u}bh\left( u\right) -2h^{2}\left( u\right) \left(
1+k_{2}b\right) +\left( b^{2}-r^{2}\left( u\right) \right) h^{\prime }\left(
u\right) ]}{h^{2}\left( u\right) }]\}dudv \\
&& \\
&&-\{b_{v}^{2}[\left( b_{u}+\dfrac{k_{1}\left( b^{2}-r^{2}\left( u\right)
\right) }{2h\left( u\right) }+k_{2}h\left( u\right) \right) ^{2} \\
&&+\dfrac{\left( k_{1}b+h^{\prime }\left( u\right) \right) [2h^{2}\left(
u\right) -2b_{u}bh\left( u\right) -2h^{2}\left( u\right) \left(
1+k_{2}b\right) +\left( b^{2}-r^{2}\left( u\right) \right) h^{\prime }\left(
u\right) ]}{h^{2}\left( u\right) }]\}dv^{2}.
\end{eqnarray*}
From this equation, we obtain the following differential equation of the loxodrome:
\begin{eqnarray*}
&&\{\left[ b_{v}\left( b_{u}+\dfrac{k_{1}\left( b^{2}-r^{2}\left( u\right)
\right) }{2h\left( u\right) }+k_{2}h\left( u\right) \right) +\frac{%
bb_{v}\left( -k_{1}b+h^{\prime }\left( u\right) \right) }{h\left( u\right) }%
\right] ^{2} \\
&&+\sinh^{2}\varphi b_{v}^{2}\left[ \left( b_{u}+\dfrac{k_{1}\left(
b^{2}-r^{2}\left( u\right) \right) }{2h\left( u\right) }+k_{2}h\left(
u\right) \right) ^{2}\right.  \\
&&\left. +\dfrac{\left( k_{1}b+h^{\prime }\left( u\right) \right)
[2h^{2}\left( u\right) -2b_{u}bh\left( u\right) -2h^{2}\left( u\right)
\left( 1+k_{2}b\right) +\left( b^{2}-r^{2}\left( u\right) \right) h^{\prime
}\left( u\right) ]}{h^{2}\left( u\right) }\right] \}\left( \dfrac{dv}{du}%
\right) ^{2} \\
&& \\
&&+\{2\cosh^{2}\varphi \left[ b_{v}\left( b_{u}+\dfrac{k_{1}\left(
b^{2}-r^{2}\left( u\right) \right) }{2h\left( u\right) }+k_{2}h\left(
u\right) \right) +\frac{bb_{v}\left( -k_{1}b+h^{\prime }\left( u\right)
\right) }{h\left( u\right) }\right]  \\
&&\times\left[ \left( b_{u}+\dfrac{k_{1}\left( b^{2}-r^{2}\left( u\right) \right)
}{2h\left( u\right) }+k_{2}h\left( u\right) \right) ^{2}\right.  \\
&&\left. +\dfrac{\left( k_{1}b+h^{\prime }\left( u\right) \right)
[2h^{2}\left( u\right) -2b_{u}bh\left( u\right) -2h^{2}\left( u\right)
\left( 1+k_{2}b\right) +\left( b^{2}-r^{2}\left( u\right) \right) h^{\prime
}\left( u\right) ]}{h^{2}\left( u\right) }\right] \}\dfrac{dv}{du} \\
&& \\
&=&-\cosh ^{2}\varphi \{\left( b_{u}+\dfrac{k_{1}\left( b^{2}-r^{2}\left(
u\right) \right) }{2h\left( u\right) }+k_{2}h\left( u\right) \right) ^{2} \\
&&+\dfrac{\left( k_{1}b+h^{\prime }\left( u\right) \right) [2h^{2}\left(
u\right) -2b_{u}bh\left( u\right) -2h^{2}\left( u\right) \left(
1+k_{2}b\right) +\left( b^{2}-r^{2}\left( u\right) \right) h^{\prime }\left(
u\right) ]}{h^{2}\left( u\right) }\}^{2}.
\end{eqnarray*}

\newpage

\section{An Example}

In this section, we give an example of the first canal surface and space-like loxodrome obtained in Section 3.

\begin{example}\rm
Let us consider the space-like spine curve $\alpha(u)=(0,0,u)$ with space-like principal normal. Taking $r(u)=u$, $m_{1}=m_{2}=1$, $\psi=\pi/3$, $u\in(0.4,2)$ and $u_{0}=1$, we have $v\in(-1.5871,1.2006)$. Moreover the arc-length of the space-like loxodrome is equal to $4.5255$. The space-like loxodrome, the space-like meridian ($v=0$) and the space-like canal surface are shown in Figure 1.

\begin{figure}[h]
\begin{center}$
\begin{array}{cc}
\includegraphics[scale=0.5]{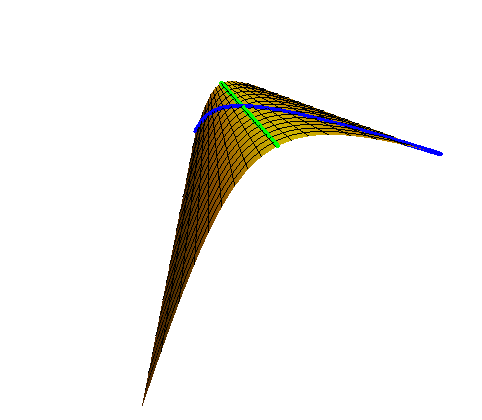}
\end{array}$
\end{center}
  \caption{Space-like loxodrome (blue), space-like meridian (green)}
\end{figure}
\end{example}

\textbf{Conflict of Interests.} The authors declare that
there is no conflict of interests regarding the publication of this paper.

\end{document}